\documentclass{amsart}

\usepackage{amssymb, amsmath, amsthm, hyperref, euscript}

\usepackage[english]{babel}

\author{Rafael Torres}

\title[Cauchy surfaces and smooth structures of $(n + 1)$-spacetimes]{Cauchy surfaces and diffeomorphism types of globally hyperbolic spacetimes}

\address{Scuola Internazionale Superiori di Studi Avanzati\\ Via Bonomea 265\\34136\\Trieste\\Italy}

\email{rtorres@sissa.it}

\dedicatory{\`a la m\'emoire de Frida}

\theoremstyle{plain}

\newtheorem{theorem}{Theorem}

\newtheorem{proposition}{Proposition}

\newtheorem{remark}{Remark}

\newtheorem{lemma}{Lemma}

\theoremstyle{definition}

\newtheorem{definition}{Definition}

\newtheorem{example}{Example}

\newcommand{\R}{\mathbb{R}}

\newcommand{\Z}{\mathbb{Z}}

\newcommand{\N}{\mathbb{N}}

\begin{document}

\maketitle

\emph{Abstract}: Chernov-Nemirovski observed that the existence of a globally hyperbolic Lorentzian metric on a $(3 + 1)$-spacetime pins down a smooth structure on the underlying 4-manifold. In this paper, we point out that the diffeomorphism type of  a globally hyperbolic $(n + 1)$-spacetime is determined by the h-cobordism class of its Cauchy surface, hence extending Chernov-Nemirovski's observation to arbitrary dimensions. 

\section{Introduction and main results}



The existence of a globally hyperbolic Lorentzian metric on a $(3 + 1)$-spacetime with closed Cauchy surface excludes all but one differentiable structure on the underlying manifold, as observed by Chernov-Nemirovski \cite{[CN]}. We point out in this note that the diffeomorphism type of a globally hyperbolic $(n + 1)$-spacetime is determined by the h-cobordism class of its closed Cauchy surface. The precise statement is as follows.

\begin{theorem} {\label{Theorem S}} Let $X$ and $X'$ be two globally hyperbolic $(n + 1)$-spacetimes with closed connected oriented Cauchy surfaces $M$ and $M'$, respectively. The following conditions are equivalent \begin{enumerate} \item $X$ is diffeomorphic to $X'$.
\item $M$ is h-cobordant to $M'$. \end{enumerate}
\end{theorem}

Theorem \ref{Theorem S} extends the aforementioned result of Chernov-Nemirovski. A decomposition theorem of Bernal-S\'anchez \cite{[BS1], [BS2]} for globally hyperbolic spacetimes is a key ingredient in its proof. Chernov-Nemirovski build upon results of Perelman \cite{[P1], [P2], [P2]} and Turaev \cite{[T]}. The crucial ingredient in the proof of Theorem \ref{Theorem S} is the s-cobordism theorem of Barden, Mazur and Stallings (see Kervaire \cite{[Ke]}, L\"uck-Kreck \cite{[LK]}).\\

There are examples of pairs of non-homeomorphic non-compact Cauchy surfaces that yield diffeomorphic $(3 + 1)$-spacetimes due to Newman-Clark \cite{[NeC]} and Chernov-Nemirovski \cite[Remark 2.3]{[CN]}. Theorem \ref{Theorem S} points out the following restrictions on the smooth structures on high-dimensional spacetimes that admit a globally hyperbolic Lorentzian metric. First, there exist infinitely many closed four-dimensional Cauchy surfaces that are homeomorphic but non-diffeomorphic whose respective globally hyperbolic $(4 + 1)$-spacetimes are diffeomorphic. In dimension at least five, the existence of a globally hyperbolic Lorentzian metric does not impose any restriction on the smooth structure of spacetime. Moreover, there are examples of homotopy equivalent spacetimes with homotopy equivalent yet non-homeomorphic Cauchy surfaces as well. A precise statement describing the phenomena is the following proposition.

\begin{proposition}{\label{Proposition S}} \begin{enumerate} \item There exist infinite sets $\{M_k: k \in \N\}$ that consist of pairwise non-diffeomorphic homeomorphic four-dimensional closed oriented Cauchy surfaces whose respective globally hyperbolic $(4 + 1)$-spacetimes are diffeomorphic to $M_1\times \R$.

\item There exist finite sets $\{M_i: i \in \{1, \ldots, k\}\}$ of cardinality $k > 1$ that consist of pairwise non-diffeomorphic homeomorphic n-dimensional closed oriented Cauchy surfaces for $n\geq 5$ whose respective globally hyperbolic $(n + 1)$-spacetimes $M_i\times \R$ are homeomorphic to $M_1\times \R$ $\forall i$, yet $M_i\times \R$ is diffeomorphic to $M_j\times \R$ if and only if $i = j$, for $j\in \{1, \ldots, k\}$.

\item There exist finite sets $\{M_i: i \in \{1, \ldots, k\}\}$ of cardinality $k > 1$ that consist of pairwise non-homeomorphic homotopy equivalent n-dimensional closed oriented Cauchy surfaces for $n\geq 5$ whose respective globally hyperbolic $(n + 1)$-spacetimes $M_i\times \R$ are homotopy equivalent to $M_1\times \R$ $\forall i$, yet $M_i\times \R$ is homeomorphic to $M_j\times \R$ if and only if $i = j$, for $j\in \{1, \ldots, k\}$.\end{enumerate}
\end{proposition}

Proposition \ref{Proposition S} is a corollary of the following results in manifold theory. The existence of manifolds that are homeomorphic yet not diffeomorphic to the n-sphere with $n\geq 7$ was proven by Milnor \cite{[KM]}. Wall \cite[Theorem]{[W]}, Hsiang-Shaneson \cite[Theorem A]{[HS]}, and Casson have shown the existence of inequivalent structures on tori of dimension at least five. There are only finitely many inequivalent smooth structures on closed manifolds of dimension greater or equal to five \cite[Theorem 1.1.8]{[GS]}. Hence, the phenomena described in Item (1) of Proposition \ref{Proposition S} can only occur in $(4 + 1)$-spacetimes (please see Example \ref{Example 4D}); we refer the reader to Gompf-Stipsicz \cite{[GS]} for details regarding the existence of inequivalent smooth structures on 4-manifolds.\\

The paper is organized as follows. In the following section, we recollect the fundamental results we are building upon to prove our main theorem. Section \ref{Section C1} contains a recollection of basic facts on globally hyperbolic spacetimes, including the fundamental result of Bernal-S\'anchez. In Section \ref{Section C2}, we present the statement of the s-cobordism theorem and its implication needed in this note. The proofs of Theorem \ref{Theorem S} and Proposition \ref{Proposition S} are given in Section \ref{Section P}.



\subsection{Acknowledgements:} We thank Vladimir Chernov and the anonymous referees for their suggestions and their careful reading of a previous version of the manuscript, which led to its substantial improvement. We thank Dionigi Benincasa, Jacobo L\'opez-Pav\'on, and Matteo Martinelli for useful discussions. 

\section{Background results}

In this section, we state the fundamental theorems regarding spacetimes and cobordisms that are needed to prove the main results of this paper.

\subsection{Decomposition of globally hyperbolic spacetimes}{\label{Section C1}} We recollect some basic notions of general relativity in this section in order to make the paper self-contained, following Hawking-Ellis' book \cite{[HE]}, and the exposition and notation in \cite[Section 1]{[CN]}. A connected manifold X equipped with a time-oriented Lorentzian metric $g$ is said to be a spacetime, and we will denote it by $(X, g)$. For any two points $x, y\in X$, we say $x\leq y$ if either $x = y$ or if there exists a future-pointing curve that connects them. A spacetime is causal if and only if there are no closed non-trivial future-pointing curves.\\

A globally hyperbolic spacetime $(X, g)$ is a causal spacetime for which the causal segments \begin{equation}I_{x, y} = \{z\in X: x\leq z\leq y\}\end{equation} are compact for every $x, y\in X$. This definition was proven by Bernal-S\'anchez \cite[Theorem 3.2]{[BS3]} to be equivalent to the classical one \cite[Section 6.6]{[HE]}. If $M\subset (X, g)$ is a subset that intersects every endless future-pointing curve exactly once, then $M$ is called a Cauchy surface. A well-known construction of globally hyperbolic spacetimes is the following. Take a complete Riemannian manifold $(M, g_r)$, and consider the metric $g_r - dt^2$ on $M\times \R$. This metric is a time-oriented Lorentzian metric for which the original Riemannian manifold is a Cauchy surface. A spacetime is globally hyperbolic if and only if it contains a Cauchy surface \cite[Section 6.6]{[HE]}. Therefore, the following well-known result is obtained.

\begin{lemma}{\label{Lemma GH}} Let $(M, g_r)$ be a complete Riemannian n-manifold. Then, \begin{equation}(M\times \R, g_r - dt^2)\end{equation} is a globally hyperbolic $(n + 1)$-spacetime.

\end{lemma}

We finish this section with the statement of the fundamental result that we build upon to prove Theorem \ref{Theorem S}.

\begin{theorem}{\label{Theorem BS}} Bernal-S\'anchez \cite{[BS1], [BS2]}. Suppose $(X, g)$ is a globally hyperbolic $(n + 1)$-spacetime. There exist a smooth n-manifold $M$ and a diffeomorphism \begin{equation}h: M \times \R \rightarrow X\end{equation} such that \begin{enumerate} \item $h(M\times \{t\})$ is a smooth spacelike Cauchy surface for all $t\in \R$, and
\item $h(\{x\}\times \R)$ is a future-pointing timelike curve for all $x\in M$. \end{enumerate}

In particular, every smooth space like Cauchy surface in $(X, g)$ is diffeomorphic to $M$.
\end{theorem}

\subsection{Diffeomorphism types of products $M\times \R$ and s-cobordisms}{\label{Section C2}} We state in this section the main topological results employed in the proof of Theorem \ref{Theorem S}. We begin by recalling the necessary notions in the following definition. The reader is referred to \cite{[Ke], [GS], [LK]} for more details on the subject.

\begin{definition} {\label{Definition S}} An (n + 1)-dimensional cobordism is a five-tuple $(W, M, f, M', f')$, where $W$ is a compact (n + 1)-manifold with boundary that decomposes as \begin{equation}\partial W = \partial_+ W \cup \partial_- W\end{equation} with orientation-preserving diffeomorphisms \begin{equation}f: \partial_-W\rightarrow -M, f':\partial_+ W \rightarrow M',\end{equation} where $-M$ denotes the manifold $M$ with the opposite orientation.\\ $\bullet$ The five-tuple is also called a cobordism $(W, M, f, M', f')$ over $M$, and $M$ is said to be cobordant to $M'$.\\
$\bullet$ Two cobordisms $(W_i, M_i, f_i, M'_i, f'_i)$ i = 1, 2 are diffeomorphic rel. $M_1$ if and only if there exists a diffeomorphism $F:W_1\rightarrow W_2$ such that $F\circ f_1 = f_1'$\\
$\bullet$ A cobordism $(W, M, f, M', f')$ over $M$ is an h-cobordism over $M$ if and only if the inclusions \begin{equation}\partial_{\pm}W\hookrightarrow W\end{equation} are homotopy equivalences.\\
$\bullet$ A cobordism $(W, M, f, M', f')$ over $M$ is trivial if and only if it is diffeomorphic to the trivial product  h-cobordism $(M\times [0, 1], M\times \{0\}, f, M\times \{1\}, f')$.
\end{definition}

The five-tuple of Definition \ref{Definition S} contains all the data necessary to employ surgery theory and the s-cobordism theorem. However, whenever there is no room for confusion and to avoid heavy notation, the maps $f$ and $f'$ in the five-tuple will be omitted and an h-cobordism will be denoted by $(W; M, M')$.

\begin{theorem}{\label{Theorem SC}} The s-cobordism theorem; Barden, Mazur, Stallings. Suppose $M$ is a closed connected smooth oriented n-manifold with $n\geq 5$ and fundamental group $\pi_1(M) = \pi$. Then, \begin{enumerate}
\item an h-cobordism $(W, M, f, M', f')$ over $M$ is trivial if and only if its Whitehead torsion $\tau(W, M) \in Wh(\pi)$ is zero;
\item for any element $\tau\in Wh(\pi)$ there exists an h-cobordism $(W, M, f, M', f')$ over $M$ such that $\tau(W, M) = \tau$;
\item the assignment of an h-cobordism $(W, M, f, M', f')$ over $M$ to its Whitehead torsion $\tau(W, M) \in Wh(\pi)$ induces a bijection between the diffeomorphism classes rel. M of h-cobordisms over $M$ to elements of $Wh(\pi)$.
\end{enumerate}
\end{theorem}

A proof of the s-cobordism theorem can be found in Kervaire \cite{[Ke]} (see \cite[Section 6]{[LK]} for more details on Whitehead torsion). 

\begin{remark}{\label{Remark L}} The s-cobordism theorem holds in dimensions $n = 1$ and $n = 2$, where all h-cobordisms are trivial. Indeed, \begin{equation}Wh(\pi_1(S^1)) = 0 = Wh(\pi_1(\Sigma_g)),\end{equation} where $\Sigma_g$ is a closed oriented surface of genus g \cite{[LK]}, Farrell-Jones \cite[Theorem 1.11]{[FJ]}. The results of Perelman \cite{[P1], [P2], [P3]} imply that the s-cobordism holds in dimension three (see \cite[Theorem 1.4]{[T]}). It is conjectured that the Whitehead group of torsion-free groups is trivial. Donaldson \cite{[D]} has shown that h-cobordant smooth oriented 4-manifolds need not be diffeomorphic.
\end{remark}

We now mention an implication of the s-cobordism theorem that is used to prove Theorem \ref{Theorem S}. 

\begin{proposition}{\label{Proposition H3}} Suppose $M$ and $M'$ are two closed connected smooth oriented n-manifolds. The following are equivalent\\
i) $M$ is h-cobordant to  $M'$.\\
ii) $M\times \R$ is diffeomorphic to $M'\times \R$.
\end{proposition}

It is straight-forward to see that Proposition \ref{Proposition H3} holds in dimension $n = 1$ given that the circle is the only closed smooth connected 1-manifold. For $n = 2$, closed oriented surface $\Sigma_g$ and $\Sigma_{g'}$ are h-cobordant if and only if they are diffeomorphic. The diffeomorphism type is determined by the fundamental group (see Remark \ref{Remark L}). The case $n = 3$ follows from work of Perelman \cite{[P1], [P2], [P3]} and Turaev \cite{[T]} as it was shown by Chernov-Nemirovski \cite[Section 3]{[CN]}. Indeed, Turaev has shown that for closed geometric oriented 3-manifolds the conditions of h-cobordant and diffeomorphic are equivalent. Perelman's solution to the Geometrization conjecture implies that all closed oriented 3-manifolds are geometric. Moreover,  it is straight forward to see that the existence of a diffeomorphism $\psi: M\rightarrow M'$ yields the existence of a diffeomorphism $\phi:M\times \R\rightarrow M'\times \R$.\\

Some words on the proof of Proposition \ref{Proposition H3} are in order. The s-cobordism theorem is invoked to prove that the implication $i)\Rightarrow ii)$ holds. Moreover, the assumption required is the existence of an h-cobordism $(W; M, f, M', f')$ for any given diffeomorphisms $f$ and $f'$. Therefore, we will omit these maps from the notation and denote an h-cobordism by $(W; M, M')$. An argument of Chernov-Nemirovski \cite[Proof Theorem B]{[CN]} is used to prove $ii)\Rightarrow i)$ and it is dimension independent.

\begin{proof} We begin by showing in two steps that the implication $i) \Rightarrow ii)$ holds. First, we show that the existence of an h-cobordism between $M$ and $M'$ implies that $M\times S^1$ and $M'\times S^1$ are diffeomorphic. In order to be able to use Theorem \ref{Theorem SC}, we need to assume that $n\geq 4$. The proof of this statement in lower dimensions has been addressed in the paragraph that follows the statement of Proposition \ref{Proposition H3}. Suppose that there exists an $(n + 1)$-dimensional h-cobordism $(V; M, M')$ where $n = dim_{\R}M = dim_{\R}M'$. We take the product $V\times S^1$, which is an $(n + 2)$-dimensional h-cobordism between $M\times S^1$ and $M'\times S^1$. The Whitehead torsion of the h-cobordism $V\times S^1$ equals $\tau(V\times S^1, M\times S^1) = 0\in Wh(\pi\times \Z)$ by \cite[Theorem 23.2]{[Co]} since the Euler characteristic of the circle is zero. Theorem \ref{Theorem SC} implies that the h-cobordism $V\times S^1$ is a trivial product h-cobordism, and that there exists a diffeomorphism \begin{equation}\varphi:M\times S^1 \rightarrow M'\times S^1.\end{equation} This concludes the first step of the proof of the claim. The second step consists of showing that the existence of the diffeomorphism $\varphi$ implies that existence of a diffeomorphism \begin{equation}{\label{Equation D}}\widetilde{\varphi}: M\times \R \rightarrow M'\times \R,\end{equation} hence concluding the proof of the implication $i)\Rightarrow ii)$. To see this, we argue as follows (cf. Hillman \cite[Lemma 6.10]{[Hi]}).. The induced isomorphism of fundamental groups \begin{equation} \varphi_{\ast}: \pi_1(M\times S^1) \rightarrow \pi_1(M'\times S^1)\end{equation} sends the fundamental group $\pi_1(M)$ to $\pi_1(M')$. Therefore, the covers \begin{equation}M\times \R\rightarrow M\times S^1\end{equation} and \begin{equation}M'\times \R\rightarrow M'\times S^1\end{equation} that correspond to the universal cover of the circle $\R\rightarrow \R/\Z = S^1$ are diffeomorphic, i.e., there exists the diffeomorphism $\widetilde{\varphi}$ of Equation \ref{Equation D}. Hence, $i)\Rightarrow ii)$, as it was claimed.

We employ an argument of Chernov-Nemirovski to see that $ii)$ implies $i)$. Suppose there exists a diffeomorphism \begin{equation}\psi: M\times \R \rightarrow M'\times \R,\end{equation}and let $T\in \R$ be sufficiently large. We use $\psi$ to build a cobordism $W$ between $M$ and $M'$. Since $M$ is compact, $\psi(M\times \{0\})$ is compact, and it is hence contained in $M'\times (-\infty, T)$. Without loss of generality, assume \begin{equation}M'\times \{T\} \subset \psi(M\times (0, +\infty)).\end{equation} We then have a smooth $(n + 1)$-manifold \begin{equation}W:= (M'\times (-\infty, T]) \cap \psi(M\times [0, + \infty))\subset M'\times \R\end{equation} whose boundary components are \begin{equation}\partial_+W:= \psi(M\times \{0\}) = M\end{equation} and \begin{equation}\partial_-W:= M'\times \{T\} = M'.\end{equation} To finish the proof, we need to check that the inclusions into $W$ of the boundary components $\partial_{\pm} W$ are homotopy equivalences. Let $r_M: M\times \R \rightarrow M\times [0, +\infty)$ and $r_{M'}: M'\times \R \rightarrow M'\times [0, +\infty)$ be the natural strong deformation retractions. A strong deformation retraction $N\times \R\rightarrow W$ is obtained by composing these maps with the diffeomorphism $\psi$ as \begin{equation} r_{M'}\circ \psi\circ r_M.\end{equation} In particular, the inclusion $W\hookrightarrow M'\times \R$ is a homotopy equivalence, as are the natural inclusions $M'\times \{T\}\hookrightarrow M'\times \R$ and $\psi(M\times \{0\}) \hookrightarrow M'\times \R$. Therefore, there exists an h-cobordism $(W; M, M')$, and the proof of Proposition \ref{Proposition H3} is now complete.

\end{proof}

\begin{remark} Kervaire \cite[Remarques]{[Ke]} provides different proofs to some statements of Proposition \ref{Proposition H3}.
\end{remark}

The following example sheds light on the avail of Proposition \ref{Proposition H3} for our purposes. 

\begin{example} \emph{Non-diffeomorphic h-cobordant manifolds, and stabilizations.} Let $L(p, q)$ be a 3-dimensional lens space whose fundamental group is the finite cyclic group $\Z/p\Z$, and denote the 2-sphere by $S^2$. Milnor has shown that, although the 5-manifolds \begin{center}$L(7, 1)\times S^2$ and $L(7, 2)\times S^2$\end{center} are h-cobordant, they are not diffeomorphic \cite{[Mi]}. In particular, this example has Whitehead group $Wh(\Z/7\Z) \neq 0$, and any h-cobordism between these two 5-manifolds is not a product h-cobordism. Proposition \ref{Proposition H3} states that the stabilizations \begin{center} $(L(7, 1)\times S^2)\times \R$ and $(L(7, 2)\times S^2)\times \R$\end{center} are however diffeomorphic.

\end{example}

\section{Globally hyperbolic $(n + 1)$-spacetimes}{\label{Section P}}

We now prove Theorem \ref{Theorem S} and Proposition \ref{Proposition S} recollecting the results presented in Section \ref{Section C1} and Section \ref{Section C2}.

\subsection{h-cobordisms of Cauchy surfaces and diffeomorphisms of spacetimes: proof of Theorem \ref{Theorem S}} Let $X$ and $X'$ be two globally hyperbolic $(n + 1)$-spacetimes with closed Cauchy surfaces $M$ and $M'$, respectively. Theorem \ref{Theorem BS} implies that there are diffeomorphisms \begin{equation} X = M\times \R,\end{equation} \begin{equation}X' = M'\times \R.\end{equation} The claim now follows from Proposition \ref{Proposition H3}.

\subsection{Myriad of Cauchy surfaces: proof of Proposition \ref{Proposition S}} Let us begin by proving Item (1). Regarding the existence of an infinite set $\{M_k: k \in \N\}$ that consists of pairwise non-diffeomorphic smooth structures on the homeomorphism type of a closed smooth oriented simply connected 4-manifold, please see Example \ref{Example 4D}; the reader is referred to \cite{[GS]} for details. Let $M_i$ be such a manifold, and $g_r$ a Riemannian metric on it. Lemma \ref{Lemma GH} implies that $(M_i\times \R, g_r - dt^2)$ is a globally hyperbolic $(4 + 1)$-spacetime for which $M_i$ is a Cauchy surface. A result of Wall states that homeomorphic simply connected smooth closed 4-manifolds are h-cobordant \cite[Theorem 9.2.13]{[GS]}. Theorem \ref{Theorem S} implies that $M_k\times \R$ are diffeomorphic for all values of $k\in \N$.

Lemma \ref{Lemma GH} and Theorem \ref{Theorem S} state that Items (2) and (3) are corollaries of the following high-dimensional results. There exists $n$-manifolds homeomorphic but not diffeomorphic to the n-torus Casson, Wall \cite[Theorem]{[W]}, Hsiang-Shaneson \cite[Theorem A]{[HS]} for $n\geq 5$. In dimensions greater than six, the statement follows from Milnor's result on the existence of manifolds homeomorphic but not diffeomorphic to the n-sphere. In particular, the h-cobordism classes of homotopy n-spheres are classified in Kervaire-Milnor \cite{[KM]}. There exist smooth manifolds homotopy equivalent to real projective n-spaces L\'opez de Medrano \cite{[Lo]}. For example, there exist four smooth 5-manifolds that are homotopy equivalent to $\mathbb{RP}^5$, and two homeomorphism classes.

\begin{example}{\label{Example 4D}} \emph{h-cobordisms and a myriad of inequivalent smooth structures in dimension four}. Dimension four is the only dimension to host closed manifolds that admit infinitely many inequivalent smooth structures. We mention now an example, and direct the interested reader towards \cite{[FS], [GS]} for further information on the subject.The hypersurface described by the homogeneous equation
\begin{equation} \{z_1^4 + z_2^4 + z_3^4 + z_4^4 = 0\}\subset \mathbb{CP}^3 \end{equation}
is diffeomorphic to the K3 surface, a well-known Calabi-Yau manifold of complex dimension two. Fintushel-Stern \cite{[FS]} have constructed an infinite set of pairwise non-diffeomorphic 4-manifolds that are homeomorphic to the K3 surface.
\end{example}

\end{document}